\documentclass[11pt]{article}
\usepackage{amssymb,amsmath,amsthm,bm,setspace}
\usepackage{pdfsync}

\newcommand{\R}{\mathbf{R}}

\renewcommand{\P}{\mathrm{P}}
\newcommand {\E}{\mathrm{E}}

\renewcommand{\d}{\text{\rm d}}
\newcommand{\e}{\text{\rm e}}
\renewcommand{\Re}{\text{\rm Re}\,}

\newcommand{\sL}{\mathcal{L}}

\newtheorem{stat}{Statement}[section]
\newtheorem{proposition}[stat]{Proposition}
\newtheorem{corollary}[stat]{Corollary}

\newtheorem{theorem}[stat]{Theorem}

\newtheorem{lemma}[stat]{Lemma}
\theoremstyle{definition}
\newtheorem{definition}[stat]{Definition}\newtheorem{remark}[stat]{Remark}
\newtheorem{example}[stat]{Example}

\numberwithin{equation}{section}

\begin{document}\onehalfspacing

\title{\bf Dynkin's isomorphism theorem
	and the stochastic heat equation%
	\thanks{%
	Research supported in part by NSF grant DMS-0704024.}}
	
\author{Nathalie Eisenbaum\\ CNRS - Universit\'e  Paris 6 \and Mohammud Foondun \\
	Loughborough University \and Davar Khoshnevisan\\ University of Utah}

\date{February 23, 2010}
\maketitle
\begin{abstract}
	Consider the stochastic heat equation
	$\partial_t u = \sL  u + \dot{W}$, where $\sL$ is
	the generator of a [Borel right] Markov process in duality. We show that
	the solution is locally mutually absolutely continuous with respect to
	a smooth perturbation of the Gaussian process that is associated,
	via Dynkin's isomorphism theorem,
	to the local times of the replica-symmetric process 
	that corresponds to $\sL$.
	In the case that $\sL$ is the generator of a L\'evy process on $\R^d$,
	our result gives a probabilistic explanation of
	the recent findings of Foondun et al.\ \cite{FKN}.\\
	
	\noindent{\it Running Title:} Dynkin's isomorphism theorem and 
	the heat equation.
\end{abstract}

\section{Introduction and main results}
The purpose of this article is to give some probabilistic
insight into the structure of the linear stochastic heat equation
\begin{equation}\left|\label{heat}\begin{split}
	&\frac{\partial}{\partial t}U(t\,,x) = (\sL U)(t\,,x) + \dot{W}(t\,,x),\\
	&U(0\,,x)=0,
\end{split}\right.\end{equation}
where $\sL$---the generator of a nice Markov process with
values on a nice space $E$---acts on the variable $x\in E$,
$t$ is strictly positive, and $\dot{W}$ is a suitable space-time
white noise on $\R_+\times E$.

A typical example is when $E=\R^d$, and $\sL$ is the $L^2$-generator
of a L\'evy process $\{X_t;\, t\ge 0\}$ on $\R^d$.
Let $X^*$ denote an independent copy of the L\'evy process $-X$ and consider
the symmetric L\'evy process $\bar X$ defined by
\begin{equation}\label{eq:barX}
	\bar X_t := X_t+X_t^*\qquad\text{for all $t\ge 0$}.
\end{equation}

It has been shown recently in \cite{FKN} that, under these conditions
on $\sL$ and $E$, \eqref{heat} has a random-field
solution $U$ if and only if $\bar X$ has local times $\{L_t^x\}_{t\ge 0,x\in\R^d}$.
Moreover, when the local times exist, many of the local features of
$x\mapsto U(t\,,x)$ are precisely the same as the corresponding festures of
$x\mapsto L_t^x$. Most notably, $x\mapsto U(t\,,x)$ is [H\"older] continuous if and only
$x\mapsto L_t^x$ is [H\"older] continuous. And the critical H\"older exponent
of $x\mapsto U(t\,,x)$ is the same as that of $x\mapsto L_t^x$.

The approach taken in \cite{FKN} is a purely analytic
one: One derives necessary and sufficient
analytic conditions for the desired local properties of $x\mapsto U(t\,,x)$ and/or
$x\mapsto L_t^x$, and checks that the analytic conditions are the same. 

The purpose of the present paper is to give an abstract probabilistic 
explanation for the many connections
that exist between the solution to \eqref{heat} and local times of the symmetrized
process $\bar X$. Our explanation does not require us to study special local properties,
and, moreover, allows us to study a much more general family of operators $\sL$ than
those that correspond to L\'evy processes. 

We close the Introduction by describing our main findings. Before we do that
we identify precisely the family operators $\sL$ with which we are concerned,
since this requires some careful development. We refer
the reader to the recent monograph by Marcus and Rosen \cite{MR},
which contains a wealth of information on Markov processes, local times,
and their deep connections to Gaussian processes. Our notation 
for Markov processes is standard and follows \cite{MR} as well.

Let $X:=\{X_t\}_{t\ge 0}$ denote a Borel right process with 
values on a locally compact, separable metric
space $E$, and let $\{P_t\}_{t\ge 0}$ denote the semigroup of $X$.
We assume that there exists a Radon measure $m$ on $E$
with respect to which $X$ has regular transition functions $p_t(x\,,y)$.

Let $L^2(m)$ denote the collection of all Borel-measurable
functions $f:E\to\R$ such that $\|f\|<\infty$, where
\begin{equation}
	\|f\|:=(f\,,f)^{1/2}\quad\text{and}\quad
	(g\,,h):=\int gh\,\d m.
\end{equation}
As usual, we define the $L^2$-generator $\sL$ and its domain
as follows: 
\begin{equation}
	\text{Dom}[\sL]:=\left\{\phi\in L^2(m):\
	\sL \phi:=\lim_{t\to 0^+}
	\frac{P_t \phi-\phi}{t}\text{ exists in $L^2(m)$}\right\}.
\end{equation}

We assume that the process $X$ has a dual process $X^*$ under $m$, so that
the adjoint  $P_t^*$ of $P_t$ is itself a Markov semigroup on $L^2(m)$. 
We emphasize that our assumptions imply
that each $P_t$ [and also $P_t^*$] is a contraction on $L^2(m)$. 
Here, and throughout, we assume also the following commutation property:
\begin{equation}\label{sym}
	P_tP_s^*=P_s^*P_t
	\qquad\text{for all $s,t\ge 0$}.
\end{equation}
This condition is met is $X$ is a L\'evy process on an abelian group, or if it is a Markov process
with symmetric transition functions.

Next, define
\begin{equation}
	\bar P_t := P_t^*P_t\qquad
	\text{for all $t\ge 0$}.
\end{equation}
A moment's thought shows that $\{\bar P_t\}_{t\ge 0}$ is
a symmetric Markovian semigroup on $L^2(m)$ simply because
of \eqref{sym}; 
$\{\bar P_t\}_{t\ge 0}$ is the \emph{replica semigroup}
associated to the process $X$, and
appears prominently in the work of Kardar \cite{Kardar},
for example.

Also, consider
the corresponding replica generator and its domain, viz.,
\begin{equation}
	\text{Dom}[\bar\sL]:=\left\{\phi\in L^2(m):\
	\bar \sL \phi:=\lim_{t\to 0^+}
	\frac{\bar P_t \phi-\phi}{t}\text{ exists in $L^2(m)$}\right\},
\end{equation}
as well as the $\alpha$-potentials $\bar U_\alpha:=
\int_0^\infty\e^{-\alpha s}\bar P_s\,\d s$ for $\alpha>0$.

Throughout, we assume that the semigroup
$\{\bar P_t\}_{t\ge 0}$ corresponds to a Borel right Markov process
$\{\bar X_t\}_{t\ge 0}$, and
$\text{Dom}[\bar\sL]$ is dense in $L^2(m)$. A simple computation
shows that $\{\bar P_t\}_{t\ge 0}$ has regular [and \emph{symmetric}]
transition functions that are denoted by $\bar p_t(x\,,y)$.

The process $\{\bar X_t\}_{t\ge 0}$ has $\alpha$-potential densities
that are described as follows: For all $\alpha>0$ and $x,y\in E$,
\begin{equation}
	\bar u_\alpha(x\,,y) =\int_0^\infty \e^{-\alpha s} \bar p_s(x\,,y)
	\,\d s.
\end{equation}
Since $\bar p_s(x\,,y)=\bar p_s(y\,,x)$,
$\{\bar X_t\}_{t\ge 0}$ is a strongly symmetric Markov process.

We are interested mainly in the case that 
$\bar u_\alpha(x\,,x)<\infty$ for all $x$ because that is precisely
the condition that guarantees that $\bar{X}$ has local times.
As we shall see [Theorem \ref{existence}], this condition is equivalent to the existence of
an a.s.-unique mild solution to \eqref{heat}. When  this condition is satisfied, we choose to
normalize the local times so that 
for all $\alpha>0$ and $x,y\in E$,
\begin{equation}
	\bar u_\alpha(x\,,y)
	=\alpha  \int_0^\infty \e^{-\alpha s} \E^x L^y_s\,\d s.
\end{equation}

In broad terms, the Dynkin isomorphism theorem \cite[Chapter 8]{MR} tells us
that many of the local properties of the local-time process $x\mapsto L^x_t$,
where $t>0$ is fixed, are the same as those of the process $\eta_\alpha$, where
$\eta_\alpha$ is a centered Gaussian process, indexed by $E$, with covariance
\begin{equation}
	\text{Cov}(\eta_\alpha(x)\,,\eta_\alpha(y))=\bar u_\alpha(x\,,y)
	\qquad\text{for all $x,y\in E$}.
\end{equation}
Here, $\alpha>0$ is a fixed but arbitrary. 

The following is the main result
of this paper.  Its proof is a combination of results proved in Sections 2 and 3.

\begin{theorem}\label{th:main} Assume that  $\bar u_\alpha$ is finite and continuous on $E\times E$
for some (equivalently for all) $\alpha>0$.  
	Let $U$ denote the unique mild solution to \eqref{heat},
	where the white noise $\dot{W}$ is chosen so that its
	control measure is $\d t\times\d m$.
	Choose and fix $\alpha>0$ and $t>0$.
	Then, there exists a space-time process $V_\alpha$
	with the following properties: 
	\begin{enumerate}
		\item For every compact set $A\subset E$, the law of
			$\{V_\alpha(t\,,x)\}_{x\in A}$ is mutually absolutely continuous with
			respect to the law of $\{U(t\,,x)\}_{x\in A}$;
		\item There exists a process $\{S_\alpha(t\,,x)\}_{x\in E}$,
			independent of $V_\alpha(t\,,\cdot)$, such that
			$S_\alpha(t\,,\cdot)$ is ``smoother than'' $V_\alpha(t\,,\cdot)$
			and $S_\alpha(t\,,\cdot)+V_\alpha(t\,,\cdot)$ has the same law
			as Dynkin's Gaussian process $\eta^\alpha$ that is associated to the
			local times of $\bar X$.
	\end{enumerate}
\end{theorem}

We will define ``smoother than'' more precisely in due time. But suffice it
to say that because $S_\alpha(t\,,\cdot)$ is ``smoother than'' $V_\alpha(t\,,\cdot)$,
many of the local properties of $S_\alpha(t\,,\cdot)$
follow from those of $V_\alpha(t\,,\cdot)$. For instance,
the following properties hold [and many more]: 
\begin{itemize}
	\item[-] If $V_\alpha(t\,,\cdot)$ is [H\"older] continuous up to a modification,
		then so is $S_\alpha(t\,,\cdot)$;
	\item[-] The critical [global/local] H\"older exponent of $S_\alpha(t\,,\cdot)$
		is at least that of $V_\alpha(t\,,\cdot)$, etc.
\end{itemize}
By mutual absolute continuity, and thanks to Dynkin's isomorphism theorem
\cite[Chapter 8]{MR}, it follows that many of the local features of
$L^\cdot_t$ and $U(t\,,\cdot)$ are shared. This explains the aforementioned
connections between \eqref{heat} and local times in the case that $\sL$
is the generator of a L\'evy process. 

We will also prove [Section 4] that when $\sL$ is
the generator of a nice L\'evy process, then we can select a $C^\infty$
version of $S_\alpha(t\,,\cdot)$. Thus, in such cases, ``smooth'' has the
usual meaning. 

Note that all the required assumptions for Theorem \ref{th:main} are satisfied   in case $\sL$ is the generator of a strongly symmetric Markov process
$X$ with finite continuous $\alpha$-potential densities ($\alpha >0$). Indeed in that case 
$\{\bar X_t\}_{t\ge 0}$ has the same law as  $\{ X_{2t}\}_{t\ge 0}$.

\section{Preliminaries}

\subsection{Markov processes}\label{subsec:Markov}
We begin by making some remarks on the underlying
Markov processes $X$ and $\bar X$. The process $X$
is chosen so that it has the following properties:
First of all, we have the identity
$(P_tf)(x) = \int_E p_t(x\,,y)f(y)\, m(\d y)$, valid
for all Borel functions $f:E\to\R_+$, $t>0$, and $x\in E$. 
And the Chapman--Kolmogorov
equation holds pointwise:
\begin{equation}
	p_{t+s}(x\,,y)=\int p_t(x\,,z)p_s(z\,,y)\, m(\d z).
\end{equation}

As was pointed out in the Introduction, a simple computation
shows that $\{\bar P_t\}_{t\ge 0}$ has regular [and \emph{symmetric}]
transition functions. In fact, they are described as follows:
for all $t>0$ and $x,y\in E$,
\begin{equation}\label{eq:pbar}
	\bar p_t(x\,,y) = \int p_t(y\,,z) p_t(x\,,z)\, m(\d z).
\end{equation}

Let us close this subsection with two technical estimates.
Here and throughout, we denote by $M(E)$
the space of finite Borel-measurable signed measures on $E$. 

\begin{lemma}\label{density} 
	If $\|P_r^*\mu\| < \infty$ for some $\mu \in M(E)$ and $r>0$, then
	$t\mapsto \|P_{t+r}^*\mu\|$ is a nonincreasing function on
	$[0\,,\infty)$.  In particular, 
	the function $t\mapsto\bar p_t(x\,,x)$ is nonincreasing on $(0\,,\infty)$
	for all $x \in E $.
\end{lemma}

\begin{proof}
	Let us choose and fix the $r>0$ as given,
	and observe that because $P_{t+r}=P_tP_r$ and $P_t$ is a contraction on $L^2(m)$,
	it follows that  $\|P_{t+r}^*\mu\|<\infty$ for all $t\ge 0$.
	
	Next, let us consider only $\mu\in\text{Dom}[\bar\sL]$, so that $\mu$
	is for the time being a function.
	Because $\bar P_t$ is a contraction on $L^2(m)$ for all $t\ge 0$, it follows that
	$\bar P_t\mu\in\text{Dom}[\bar\sL]$ for all $t\ge 0$, and therefore
	\begin{equation}\label{23}\begin{split}
		\frac{\d}{\d s}\left\|
			P_s^*\mu\right\|^2 =
			\frac{\d}{\d s}( P_s^*\mu\,,P_s^*\mu)  &=
			\frac{\d}{\d s}(\bar P_s\mu\,,\mu)
			=(\bar\sL \bar P_s\mu\,,\mu) \\
		&\hskip.8in=(\bar\sL P_s\mu\,,P_s\mu),
	\end{split}\end{equation}
	where $\d/\d s$ denotes the right derivative at $s$.
	It is well known that $\bar{\sL}$ is a negative-definite operator. 
	That is,
	\begin{equation}\label{negdef}
		(\bar\sL\phi\,,\phi)\le 0
		\quad\text{for all $\phi\in\text{Dom}[\bar\sL]$.}
	\end{equation}
	Indeed, because every $P_t$ is a contraction on $L^2(m)$, it follows that
	$(\bar P_t\phi\,,\phi)=\|P_t\phi\| ^2\le
	(\phi\,,\phi)$. Take the difference, divide by $t$, and 
	then let $t\downarrow 0$
	to deduce \eqref{negdef}. In turn, \eqref{negdef} and \eqref{23} 
	together imply that
	for all $\mu\in\text{Dom}[\bar\sL]$,
	\begin{equation}\label{eq:mono}
		\|P_{s+t}^*\mu\|^2 \le  \|P_s^*\mu\|^2
		\qquad\text{for all $s,t\ge 0$}.
	\end{equation}
	Since every $P_t^*$ is a contraction on $L^2(m)$, the assumed
	density of the domain of $\bar\sL$
	implies that \eqref{eq:mono} continues to hold for all $\mu\in L^2(m)$ and
	$s,t>0$.
	
	Now let $\mu$ be a finite signed Borel measure on $E$ such that 
	$\|P_r^*\mu\| < \infty$; we can apply \eqref{eq:mono},
	with $\mu$ replaced by $P_r^*\mu\in L^2(m)$, and this leads us to the following
	inequality:
	\begin{equation}
		\|P_{s+t}^*P_r^*\mu\|  \le \|P_s^*P_r^*\mu\|
		\quad\text{for all $s,t>0$}.
	\end{equation}
          Since $P_u^*P_v^*\mu=P_{u+v}^*\mu$, the preceding shows that
	\eqref{eq:mono} holds for all $s \geq r$ and all $t\geq 0$.  
	
	In order to conclude, we choose $\mu := \delta_x$. In that case, the
	Chapman--Kolmogorov equations imply that
	$ \|P_r^*\mu\| = \bar p_r(x\,,x)^{1/2} < \infty$ for all $r>0$;
	therefore, \eqref{eq:mono} implies the announced result.
	\end{proof}

 Let us end this subsection by introducing an estimate
 on $\alpha$-potentials. 
 
\begin{lemma}\label{lem:green} Assume that  $\bar u_\alpha$ is finite on $E\times E$ for $\alpha > 0$. Then
for all $x,y\in E$ and $\alpha>0$,
	\begin{equation}
		\bar u_\alpha(x\,,y) \le c_\alpha \bar u_1(y\,,y),
	\end{equation}
	where $c_\alpha=\e(\alpha+2\alpha^{-1})$.
\end{lemma}

\begin{proof}
	We begin by proving that for all $x,y\in E$ and $t\ge 1$,
	\begin{equation}\label{eq:goal:LT}
		\E^x  L^y_t\le 2t \E^yL^y_1.
	\end{equation}
	In order to prove this we recall that if $\theta$ denotes the
	shifts on the path of the underlying Markov process, then
	$L^y_{t+s}=L^y_t+L^y_s\circ\theta_t$, $\P^x$-almost surely.
	Therefore,
	\begin{equation}
		\E^x L^y_{t+s} = \E^x L^y_t +
		\E^x \E^{X_t} L^y_s.
	\end{equation}
	We can apply the strong Markov property to the first hitting time
	of $y$ to find that $\E^x L_v^y\le\E^y L^y_v$ for
	all $x,y\in E$ and $v\ge 0$. Consequently,
	\begin{equation}\label{eq:LL}
		\E^x L^y_{t+s} \le \E^y L^y_t + \E^y L^y_s
		\quad\text{for all $x,y\in E$ and $s,t\ge 0$};
	\end{equation}
	and therefore,  $\E^x L_n^y\le n\E^yL^y_1$ for all integers
	$n\ge 1$ and $x,y\in E$. If $t\ge 1$, then we can find an integer
	$n\ge 2$ such that $n-1\le t<n$, whence
	\begin{equation}
		\E^x L^y_t \le \E^x L^y_n \le n\E^y L^y_1\le (t+1)\E^y L^y_1
		\le 2t \E^y L^y_1.
	\end{equation}
	This establishes \eqref{eq:goal:LT} for $t\ge 1$.
	
	Next, we note that since $t\mapsto L_t^y$ is nondecreasing 
	a.s.\ [$\P^a$ for all $a\in E$],
	\begin{equation}\label{eq:step1}\begin{split}
		\bar u_\alpha(x\,,y)
			&= \int_0^\alpha \e^{-r} \E^x L_{r/\alpha}^y\,\d r
			+\int_\alpha^\infty \e^{-r} \E^x L^y_{r/\alpha}\,\d r\\
		&\le \alpha \E^x L^y_1+\frac{2}{\alpha}\E^yL^y_1\cdot\int_\alpha^\infty
			r\e^{-r}\,\d r\\
		&\le (\alpha + 2\alpha^{-1})\E^y L^y_1,
	\end{split}\end{equation}
	thanks to \eqref{eq:goal:LT} and the strong
	Markov property. On the other hand,
	\begin{equation}
		\bar u_1(y\,,y) = \int_0^\infty \e^{-r}\E^y L^y_r\,\d r
		\ge \E^y L^y_1 \cdot \int_1^\infty \e^{-r}\,\d r.
	\end{equation}
	The lemma follows from the preceding and \eqref{eq:step1}.
\end{proof}

\subsection{Gaussian random fields}
Suppose $\{G(x)\}_{x\in E}$ is a centered Gaussian process
that is continuous in $L^2(\P)$. The latter means that
$\E(|G(x)-G(y)|^2)\to 0$ as $x\to y$. It follows that $G$
has a separable version for which we can define
\begin{equation}
	G(\mu):=\int_E G\,\,\d\mu,
\end{equation}
for all $\mu$ in the space $M(E)$ of
finite Borel-measurable signed measures on $E$.
Moreover, $\{G(\mu)\}_{\mu\in M(E)}$ is a Gaussian
random field with mean process zero and
\begin{equation}
	\text{Cov}( G(\mu)\,,G(\nu) )= \iint \E[G(x)G(y)]\,\mu(\d x)
	\,\nu(\d y).
\end{equation}

\begin{definition}
	Let $G$ and $ G_* $ denote two $L^2(\P)$-continuous
	Gaussian processes indexed by $E$. We say that
	$G$ is \emph{smoother than} $ G_* $ if
	there exists a finite constant $c$ such that 
	\begin{equation}
		\E\left(\left| G(x)-G(y) \right|^2\right) \le
		c\E\left(\left|  G_* (x)- G_* (y)\right|^2\right)
		\qquad\text{for all $x,y\in E$}.
	\end{equation}
	We say that $G$ is \emph{as smooth as} $G^*$ when
	$G$ is smoother than $G^*$ and $G^*$ is smoother than $G$.
\end{definition}

It is easy to deduce from general theory \cite[Chapters 5--7]{MR} that 
if $G$ is smoother than $ G_* $, then 
the continuity of $x\mapsto G_* (x)$
implies the continuity of $x\mapsto G(x)$. Similar
remarks can be made about H\"older continuity and
existence of nontrivial $p$-variations, in case the
latter properties hold and/or make sense.

\section{The heat and cable equations}

Let $\dot{W}:=\{\dot{W}(t\,,x);\, t\ge 0,\, x\in E\}$ denote
white noise on $\R_+\times E$ with control measure
$\d t\times \d m$, and
consider the stochastic heat equation \eqref{heat},
where $t>0$ and $x\in E$.
Also, consider the stochastic cable equation with parameter
$\alpha>0$:
\begin{equation}\label{cable}\left|\begin{split}
	&\frac{\partial}{\partial t} V_\alpha(t\,,x) = (\sL V_\alpha)(t\,,x)
		- \frac{\alpha}{2} V_\alpha(t\,,x) + \dot{W}(t\,,x),\\
	&V_\alpha(0\,,x)=0,
\end{split}\right.\end{equation}
where $t>0$ and $x\in E$.

First we establish the existence of a mild solution to \eqref{cable}
and to \eqref{heat} [Theorem \ref{existence}]. Since these are
linear equations, issues of uniqueness do not arise.

Then, we return to the first goal of this section and prove that all local properties
of $U(t\,,\cdot)$ and $V_\alpha(t\,,\cdot)$ are the same;
see Proposition \ref{pr:girsanov}. We do this in two steps.
First, we study the case that $\alpha$
is sufficiently small. In that case, we follow a change
of measure method that is completely analogous to Proposition 
1.6 of Nualart and Pardoux \cite{NP}; see also 
Dalang and Nualart \cite[Theorem 5.2]{DalangNualart}.
In a second step, we bootstrap from small values of $\alpha$ to
large values of $\alpha$ by an argument that
might have other applications as well.

Next, we prove that $V_\alpha(t\,,\cdot)$ is 
equal to a smooth perturbation of the associated Gaussian
process that arises in the Dynkin isomorphism theorem;
see Proposition \ref{pr:bd2}. The key estimate
is a peculiar inequality [Proposition \ref{pr:bd2}] that
is nontrivial even when $X$ is Brownian motion. 

As a consequence of all this, and
thanks to the Dynkin isomorphism theorem,
many of the local properties of the solution to \eqref{heat}
are the same as those of the local times of the process $\bar X$.
We refer the reader to Chapter 8 of the book by Marcus and Rosen
\cite{MR} for details on Dynkin's isomorphism theorem and its
applications to the analysis of local properties of local times.

Let us concentrate first on the cable equation.

The weak solution to the Kolmogorov equation
\begin{equation}\left|\begin{split}
	&\frac{\partial}{\partial t} f(t\,,x\,,y) = 
		(\sL_y f)(t\,,x\,,y)-\frac{\alpha}{2}f(t\,,x\,,y),\\
	&f(0\,,x\,,\cdot)=\delta_x,\\
\end{split}\right.\end{equation}
is the function $f(t\,,x\,,y):=\e^{-\alpha t/2}p_t(x\,,y)$. Therefore,
we can use the theory of Walsh
\cite[Chapter 3]{Walsh}, and write the solution to
\eqref{cable} as
\begin{equation}\label{eq:rep1}
	V_\alpha(t\,,x) = \int_0^t \int_E \e^{-\alpha(t-s)/2}p_{t-s}(x\,,y)
	\,W(\d y\,\d s).
\end{equation}
This is a well-defined Gaussian process if and only if the
stochastic integral has two finite moments. But then,
Wiener's isometry tells us that 
\begin{equation}\label{eq:EV2}\begin{split}
	\E\left(\left| V_\alpha(t\,,x)\right|^2\right)
		&=\int_0^t\e^{-\alpha s}\d s\int_E m(\d y)\
		\left| p_s(x\,,y)\right|^2\\
	&=\int_0^t\e^{-\alpha s} \bar p_s(x\,,x)\,\d s.\\
\end{split}\end{equation}

Similarly, the stochastic heat equation \eqref{heat} 
has the following solution:
\begin{equation}
	U(t\,,x) = \int_0^t \int_E 
	p_{t-s}(x\,,y)\, W(\d y\,\d s),
\end{equation}
which is a well-defined Gaussian process if and only if its second moment is finite. Note that
\begin{equation}\label{moment}
\E\left(\left| U(t\,,x)\right|^2\right) = \int_0^t\bar p_s(x\,,x)\,\d s.
\end{equation}

\begin{theorem}\label{existence} 
	\begin{enumerate}
	\item The stochastic cable equation (\ref{cable}) has an a.s.-unique mild solution if and only if 
		the Markov process  $\bar{X}$ has local times. 
	\item The stochastic heat equation (\ref{heat}) has an a.s.-unique mild solution if and only if 
		the Markov process  $\bar{X}$ has local times. 
\end{enumerate}
\end{theorem}

\begin{proof} 
	According to \cite[Lemma 3.5]{FKN}, the following holds 
	for every nonincreasing  measurable function $g:\R_+\to\R_+$, 
	and $t, \lambda > 0$: 
	\begin{equation}\label{real}
		(1 - \e^{-2 t/\lambda}) \int_0^{\infty}\e^{-2s/\lambda} g(s)\,\d s 
		\leq \int_0^t g(s)\,\d s \leq  \e^{2 t /\lambda}\int_0^{\infty}
		\e^{-2s/\lambda} g(s)\,\d s.
	\end{equation}
	We apply Lemma \ref{density} and \eqref{real}, 
	with $\lambda:=2/\alpha$ and $g(s) := \e^{-\alpha s} \bar p_s(x\,,x)$
	to find that
	\begin{equation}
		(1 - \e^{-t \alpha}) \bar u_{2\alpha}(x\,,x) \leq  
		\E\left(\left| V_\alpha(t\,,x)\right|^2\right)
		\leq \e^{t\alpha} \bar u_{2\alpha}(x\,,x)
	\end{equation}
	Similarly, we apply \eqref{real} with $\lambda:=1/\alpha$
	and $g(s):=\bar p_s(x\,,x)$ to obtain
	\begin{equation}
		(1 - \e^{-2 t \alpha}) \bar u_{2\alpha}(x\,,x) \leq  
		\E\left(\left| U(t\,,x)\right|^2\right) \leq \e^{2 t \alpha} \bar u_{2\alpha}(x\,,x).
	\end{equation}
	This proves that the processes $V_\alpha$ and $U$ are well defined
	if and only if $u_{2\alpha}(x\,,x)<\infty$ for all $x\in E$
	[and some---hence all---$\alpha>0$].
	And the latter is equivalent to the existence of local times;
	see Blumenthal and Getoor \cite[Theorem 3.13 and (3.15), pp.\ 216--217]{BG}.
 \end{proof}

From now on,  we assume that all the $\alpha$-potentials are finite and continuous on $E\times E$. 
Hence the stochastic cable equation
\eqref{cable} has an a.s.-unique mild solution such that
$\sup_{t\ge 0}\sup_{x\in A}\E(|V_\alpha(t\,,x)|^2)<\infty$
for every compact set $A\subset E$.
It follows easily from this discussion that $x\mapsto V_\alpha(t\,,x)$ is
continuous in $L^2(\P)$, and hence in probability as well. 

Since $x\mapsto V_\alpha(t\,,x)$ is continuous in probability,
Doob's separability theorem implies that
$V_\alpha(t\,,\mu):=\int_E V_\alpha(t\,,x)\,\mu(\d x)$ is well-defined for all
finite signed Borel measures $\mu$ on $E$. We will be 
particularly interested in the two examples,
$\mu:=\delta_a$ and $\mu:=\delta_a-\delta_b$ for fixed $a,b\in E$.
In those two cases, $V_\alpha(t\,,\mu)=V_\alpha(t\,,a)$ and
$V_\alpha(t\,,\mu)=V_\alpha(t\,,a)-V_\alpha(t\,,b)$,
respectively. 
One can prove quite easily the following: For all finite signed Borel
measures $\mu$ on $E$,
\begin{equation}
	V_\alpha(t\,,\mu) = \int_0^t\int_E \e^{-\alpha(t-s)/2}
	(P_{t-s}^*\mu)(y)\, W(\d y\,\d s),
\end{equation}
provided that $(|\mu|\,,\bar U_\alpha|\mu|)<\infty$.
This can be derived by showing that the second moment of 
the difference of the two quantities is zero. It also follows from
the stochastic Fubini theorem of Walsh
\cite[Theorem 2.6, p. 296]{Walsh}.

 Using the same methods
as before, one shows that $U(t\,,\mu)$ is well defined if and only if
$\mu\in M(E)$ satisfies $(|\mu|\,,\bar U_\alpha|\mu|)<\infty$.
The following result shows that each $V_\alpha(t\,,\cdot)$
is as smooth as $U(t\,,\cdot)$. A much better result
will be proved subsequently.

\begin{lemma}\label{lem:UV}
	For all $t,\alpha>0$ and $\mu\in M(E)$ such that
	$(|\mu|\,,\bar U_\alpha|\mu|)<\infty$,
	\begin{equation}
		\E\left( \left| V_\alpha(t\,,\mu)\right|^2\right)\le
		\E\left( \left| U(t\,,\mu)\right|^2\right)\le
		3\e^{\alpha t}
		\E\left( \left| V_\alpha(t\,,\mu)\right|^2\right).
	\end{equation}
	Thus, $\{V_\alpha(t\,,x);\, x\in E\}$ is as smooth as
	$\{U(t\,,x);\, x\in E\}$ for every
	$t,\alpha>0$.
\end{lemma}

\begin{proof}
	The first inequality follows from a direct computation.
	For the second bound we note that if $0\le s\le t$, then
	$(1-\e^{-\alpha s/2})^2 \le \left(\e^{\alpha t/2}-1\right)^2 \e^{-\alpha s}
	\le \e^{\alpha (t-s)}$.
	Therefore,
	\begin{align}\nonumber
		\E\left(\left| U(t\,,\mu)-V_\alpha(t\,,\mu)\right|^2\right)
			&=\int_0^t\int_E \left(1-
			\e^{-\alpha s}\right)^2
			\left|(P_s^*\mu)(y)\right|^2\,m(\d y)\,\d s\\
		&\le \e^{\alpha t}\int_0^t \e^{-\alpha s}
			\| P_s^*\mu\| ^2\,\d s\\\nonumber
		&=\e^{\alpha t}
			\E\left(\left| V_\alpha(t\,,\mu)\right|^2\right).
	\end{align}
	Because
	\begin{equation}\begin{split}
		\E\left[ U(t\,,\mu)V_\alpha(t\,,\mu)\right] &=\int_0^t
			\e^{-\alpha s}\left\| P_s^*\mu\right\| ^2\,\d s\\
		&\le \e^{\alpha t}\int_0^t \e^{-\alpha s}\left\|
			P_s^*\mu\right\| ^2\,\d s\\
		&=\e^{\alpha t}\E\left(\left| V_\alpha(t\,,\mu)\right|^2\right),
	\end{split}\end{equation}
	the second inequality of the lemma follows.
\end{proof}

We propose to prove Proposition \ref{pr:girsanov} which is a better version of Lemma
\ref{lem:UV}. But
first recall that laws of two real-valued
random fields $\{A_v\}_{v\in \Gamma}$
and $\{B_v\}_{v\in \Gamma}$ are said to be mutually
absolutely continuous if there exists an almost surely
strictly-positive mean-one random variable
$D$ such that for every $v_1,\ldots,v_n\in\Gamma$
and  Borel sets $\Sigma_1,\ldots,\Sigma_n\in\R$,
\begin{equation}
	\P\left(\bigcap_{j=1}^n\left\{A_{v_j}\in \Sigma_j\right\}\right)
	= \E\left(D\,;\,\bigcap_{j=1}^n\left\{B_{v_j}\in \Sigma_j\right\}\right).
\end{equation}

\begin{proposition}\label{pr:girsanov}
	Choose and fix $T>0$ and a compact set $A\subset E$.
	Then, then for any $\alpha>0$, the law of the random field
	$\{V_\alpha(t\,,x)\}_{t\in[0,T],x\in A}$ is mutually absolutely continuous
	with respect to the law of the random field $\{U(t\,,x)\}_{t\in[0,T],x\in A}$.
\end{proposition}

\begin{proof}
	Throughout this proof define for all $t\in[0\,,T]$,
	\begin{equation}
		Q_\alpha(t\,,A):=\int_0^t\d s\int_Am(\d y)\
		|V_\alpha(x\,,y)|^2.
	\end{equation}
	Minkowski's inequality implies that for all integers $k\ge 0$,
	\begin{equation}
		\left\| Q_\alpha(T\,,A)\right\|_{L^k(\P)}
		\le \int_0^T\d s \int_A m(\d y) \left\| 
		V_\alpha(s\,,y)\right\|_{L^{2k}(\P)}^2.
	\end{equation}
	Because of \eqref{eq:rep1}, each $V_\alpha(s\,,y)$ is a centered Gaussian 
	random variable. Therefore,
	\begin{equation}
		\| V_\alpha(s\,,y)\|_{L^{2k}(\P)}^2= \|V_\alpha(s\,,y)\|_{L^2(\P)}^2\cdot
		\|Z\|_{L^{2k}(\P)}^2,
	\end{equation}
	where $Z$ is a standard-normal random variable. Thanks to \eqref{eq:EV2},
	\begin{equation}
		\E\left(\left| Q_\alpha(T\,,A)\right|^k\right)
		\le C_\alpha^k\cdot \E(Z^{2k}),
	\end{equation}
	where $C_\alpha:= C_\alpha (T\,,A) := Tm(A)\sup_{x\in A}\bar u_\alpha(x\,,x)$.
	Consequently, for any positive integer $l$,
	\begin{equation}
		\E\left[\exp\left(\frac{\alpha^2}{8l^2}Q_\alpha(T\,,A) \right)\right]
		\le\E\left[\exp\left( \frac{\alpha^2 C_\alpha Z^2
		}{8l^2}\right)\right],
	\end{equation}
	and this is finite if and only if $\alpha^2 C_\alpha/l^2<4$.
	By continuity, $\bar u_1$ is bounded uniformly on $A\times A$. Therefore,
	Lemma \ref{lem:green} tells us that there exist a large $l$ such that  $\alpha^2C_\alpha/l^2<4$.
	Hence for such an $l$,
	\begin{equation}
		\E\left[\exp\left(\frac{\alpha^2}{8l^2}Q_\alpha(T\,,A)\right)\right]<\infty.
	\end{equation}
	Consequently, we can apply the criteria of Novikov and Kazamaki 
	(see Revuz and Yor \cite[pp.\ 307--308]{RY})
	to conclude that
	\begin{equation}
		\exp\left( \frac{\alpha}{2l}\int_0^t \int_A V_\alpha(s\,,y)\, W(\d y\,\d s)
		-\frac{\alpha^2}{8l^2}Q_\alpha(t\,,A)\right)
	\end{equation}
	defines a mean-one martingale indexed by $t\in[0\,,T]$. 
	Define
	\begin{equation}
		\dot{W}^{(1)}(t\,,x):= \dot{W}(t\,,x) - \frac{\alpha}{2l} V_\alpha(t\,,x)
		\qquad\text{for $t\in[0,T]$ and $x\in A$.}
	\end{equation}
	Recall that $\P$ denotes the measure under which $\{\dot{W}(t\,,x)\}_{t\in[0,T],x\in A}$ 
	is a white noise. The preceding and 
	Girsanov's theorem together imply that 
	$\{\dot{W}^{(1)}(t\,,x)\}_{t\in[0,T],x\in A}$
	is a white noise under a different probability measure $\P_1$ which is 
	mutually absolutely continuous with respect to $\P$;
	see Da Prato and Zabczyk \cite[p.\ 290]{DPZ}. 
	Next, we define iteratively,
	\begin{equation}
		\dot{W}^{(n+1)}(t\,,x):=\dot{W}^{(n)}(t\,,x)-
		\frac{\alpha }{l}V_\alpha(t\,,x)\quad \text{for} \quad n=1,\ldots,l-1.
	\end{equation}
	A second application of  Girsanov's theorem allows us to conclude that 
	$\{\dot{W}^{(2)}(t\,,x)\}_{t\in[0,T],x\in A}$ is a white noise under a 
	certain probability measure $\P_2$ which is 
	mutually absolutely continuous with respect to $\P_1$. 
	
	In fact, the very same argument 
	implies existence of a finite sequence of measures $\{\P_n\}_{n=1}^{l}$ 
	such that $\{\dot{W}^{(n)}(t\,,x)\}_{t\in[0,T],x\in A}$ is a white noise 
	under $\P_n$ for $1\le n\le l$.  We can now conclude that  
	$\dot{W}^{(l)}(t\,,x)= \dot{W}(t\,,x) - (\alpha/{2}) V_\alpha(t\,,x)$ 
	defines a white noise [indexed by $t\in[0\,,T]$
	and $x\in A$] under the  measure $\P_{l}$, and that $\P_l$ 
	is mutually absolutely continuous 
	with respect to $\P$. The latter fact follows from the transitivity property of 
	absolute continuity of measures; this is the property that asserts that 
	whenever $\mathrm{Q}_1$ and $\mathrm{Q}_2$ are mutually
	absolutely continuous probability measures, and $\mathrm{Q}_2$ and $\mathrm{Q}_3$
	are mutually absolutely continuous probability measures, then so are $\mathrm{Q}_1$
	and $\mathrm{Q}_3$.
	
	The result follows from the strong existence of solutions to \eqref{heat}
	and \eqref{cable}.
\end{proof}

Consider the Gaussian random field
\begin{equation}
	S_\alpha(t\,,x) := \int_t^\infty\int_E \e^{-\alpha s/2}
	p_s(x\,,y)\,W(\d y\,\d s).
\end{equation}
One verifies, just as one does for $V_\alpha(t\,,\cdot)$, that for
all finite Borel signed measures $\mu$ on $E$,
\begin{equation}
	S_\alpha(t\,,\mu) := \int_t^\infty\int_E\e^{-\alpha s/2}
	(P_s^*\mu)(y)\, W(\d y\,\d s),
\end{equation}
indexed by $\mu\in L^2(m)$ and $t>0$.
Elementary properties of the processes $S_\alpha$ and $V_\alpha$ show
that they are independent mean-zero Gaussian processes.
Consider the Gaussian random field
\begin{equation}\label{eq:VZ}
	\eta_\alpha(t\,,\cdot) := V_\alpha(t\,,\cdot)+S_\alpha(t\,,\cdot).
\end{equation}
Because
$\E [\eta_\alpha(t\,,\mu)\eta_\alpha(t\,,\nu) ]
= \E [ V_\alpha(t\,,\mu)V_\alpha(t\,,\nu) ] + 
\E [ S_\alpha(t\,,\mu)S_\alpha(t\,,\nu) ]$,
a direct computation shows that for all $t>0$
and finite Borel measures $\mu$ and $\nu$ on $E$,
\begin{equation}\label{eq:norm:eta}\begin{split}
	\text{Cov}\left(\eta_\alpha(t\,,\mu)\,,\eta_\alpha(t\,,\nu)\right)
		&=\int_0^\infty \e^{-\alpha s} \left( P_s^*\mu
		\,,P_s^*\nu\right) \,\d s\\
	&=\int_0^\infty \e^{-\alpha s} \left( \mu
		\,,\bar P_s\nu\right) \,\d s\\
	&=(\mu\,,\bar U_\alpha \nu).
\end{split}\end{equation}

In other words, the law of $\eta_\alpha(t\,,\cdot)$ does not depend on $t>0$, and
\begin{equation}
	\text{Cov}(\eta_\alpha(t\,,x)\,,\eta_\alpha(t\,,z))=\bar u_\alpha(x\,,z).
\end{equation}
Thus, $\eta_\alpha(t\,,\cdot)$ is precisely the associated Gaussian process that
arises in Dynkin's isomorphism theorem \cite[Chapter 8]{MR}.

It is easy to see that since the law of $\eta_\alpha(t\,,\cdot)$ is independent of
$t>0$, $\eta_\alpha(t\,,\cdot)$ is the [weak]
steady-state solution to \eqref{cable}, in the sense that 
$S_\alpha(t\,,x)\to 0$ in $L^2(\P)$ as $t\to\infty$ for all $x\in E$;
this follows directly from the definition of $S_\alpha$.

Our next result implies that many of the local
regularity properties of $V_\alpha(t\,,\cdot)$ and 
$\eta_\alpha(t\,,\cdot)=V_\alpha(t\,,\cdot) +S_\alpha(t\,,\cdot)$
are shared.

\begin{proposition}\label{pr:bd2}
	For every fixed $t,\alpha>0$,
	$S_\alpha(t\,,\cdot)$ is smoother than
	$V_\alpha(t\,,\cdot)$. In fact, for all 
	$\mu\in M(E)$ such that $(|\mu|\,,\bar U_\alpha |\mu|)<\infty$,
	\begin{equation}
		\E\left(\left| S_\alpha(t\,,\mu)\right|^2\right) \le
		\left[\frac{1}{\e^{t\alpha}-1}\right]\cdot
		\E\left(\left| V_\alpha(t\,,\mu)\right|^2\right).
	\end{equation}
\end{proposition}

\begin{proof}
	Suppose $\mu$ is a finite
	signed Borel measure on $E$ with $(|\mu|\,,\bar U_\alpha |\mu|)<\infty$.
	Then, $\int_0^\infty \e^{-\alpha s} \|P_s^*\mu\| ^2\,\d s
	=(\mu\,,\bar U_\alpha \mu)<\infty$, and hence
	$\|P_r^*\mu\|$ is finite for almost all $r>0$.
	
	We first note the following: 
	\begin{equation}\label{eq:328}
		\int_t^\infty \e^{-\alpha s}\left\|
		P_s^*\mu\right\|^2 \, \d s
		=\sum_{n=1}^\infty \e^{-n\alpha t}\cdot\int_0^t
		\e^{-\alpha s}\left\|
		P_{s+nt}^*\mu\right\|^2 \,\d s.
	\end{equation}
	
	But thanks to Lemma \ref{density}, we have 
	\begin{equation}
		\|P_{s+nt}^*\mu\|  \le \|P_s^*\mu\|
		\quad\text{for all $s,t>0$ and $n\ge 1$},
	\end{equation}		
	which together with \eqref{eq:328} implies that
	\begin{equation}\label{eq:strange}
		\int_t^\infty \e^{-\alpha s}\left\|
		P_s^*\mu\right\|^2 \, \d s 
		\le \left[\frac{\e^{-t\alpha}}{1-\e^{-t\alpha}}\right]\cdot
		\int_0^t\e^{-\alpha s}\left\|
		P_s^*\mu\right\|^2 \,\d s.
	\end{equation}
	This is another way of stating the Proposition.
\end{proof}

We mention, in passing, a nontrivial consequence of Proposition \ref{pr:bd2}:
Consider the special case that $\mu:=\delta_a-\delta_b$ for fixed
$a,b\in E$. In that case,
$(|\mu|\,,\bar U_\alpha|\mu|)$ is finite. Indeed,
$(|\mu|\,,\bar U_\alpha|\mu|)=\bar u_\alpha(a\,,a)
+\bar u_\alpha(b\,,b) -2\bar u_\alpha(a\,,b)$, from which it
follows that
$\|P^*_s\mu\|^2=\bar p_s(a\,,a) + \bar p_s(b\,,b)-2\bar p_s(a\,,b)$.
Therefore, time reversal and
Proposition \ref{pr:bd2}---specifically in the form given by
\eqref{eq:strange}---together assert the following somewhat unusual 
inequality.

\begin{corollary}\label{co:bd2}
	Let $S(\alpha)$ denote an independent exponentially-distributed
	random variable with mean $1/\alpha$. Then, for all
	$t,\alpha>0$ and $a,b\in E$,
	\begin{equation}
		\E^a\left[\left. L^a_{S(\alpha)}-L^b_{S(\alpha)}\ \right|\, S(\alpha)\ge t\right]
		\le 
		\E^a\left[\left. L^a_{S(\alpha)}-L^b_{S(\alpha)}\ \right|\, S(\alpha)<t\right].
	\end{equation}
\end{corollary}

This appears to be novel even when $X$ is linear Brownian motion.

\section{L\'evy processes}

Next we study the special case
that $\sL$ is the generator of a L\'evy process
$X$ on $\R$. Recall that $X$ is best understood via
its characteristic exponent $\Psi$ \cite{Bertoin}; we
normalize that exponent as follows:
$\E\exp(i\xi X_t)=\exp(-t\Psi(\xi))$. 

Let $m$ denote the Lebesgue measure on $E:=\R$; then $X$
is in duality with $-X$ under $m$, and the replica semigroup
$\{P^*_t\}_{t\ge 0}$ is the semigroup associated to the
L\'evy process $\bar{X}$ defined by \eqref{eq:barX}.

We know from the general theory of Dalang \cite{Dalang}
that \eqref{heat} has a random-field solution if and only if
\begin{equation}\label{eq:dalang}
	\int_{-\infty}^\infty \frac{\d\xi}{\alpha+2\Re\Psi(\xi)}<\infty,
\end{equation}
for one, and hence all, $\alpha>0$. This condition implies the existence
of jointly-continuous
transition functions for both $X$ and $\bar X$ \cite[Lemma 8.1]{FKN}, and
\begin{equation}
	\bar u_\alpha(x\,,y) = \frac{1}{\pi}\int_0^\infty
	\frac{\cos(\xi(x-y))}{\alpha+2\Re\Psi(\xi)}\,\d\xi.
\end{equation}
In particular, $\bar u_\alpha$ is continuous on $\R\times\R$.
Finally, one checks that the domain of $\bar\sL$ is precisely 
the collection of all $f\in L^2(m)$ such that $\Re\Psi\cdot|\hat f|^2
\in L^1(\R)$. The well-known fact that $\Re\Psi(\xi)=O(\xi^2)$
as $|\xi|\to\infty$ tells us that
all rapidly-decreasing test functions are in $\text{Dom}[\bar\sL]$,
and therefore $\text{Dom}[\bar\sL]$ is dense in $L^2(m)$.
Thus, all the conditions of Theorem \ref{th:main} are verified in
this case. 

Choose and fix some $t>0$.
According to Proposition \ref{pr:bd2}, and thanks to the general
theory of Gaussian processes, the process $S_\alpha(t\,,\cdot)$
is at least as smooth as the process $V_\alpha(t\,,\cdot)$; the
latter solves the stochastic cable equation. Next we prove that
under a mild condition on $\Psi$, $S_\alpha(t\,,\cdot)$ is in fact
extremely smooth. 

\begin{proposition}\label{pr:smooth}
	Suppose, in addition to \eqref{eq:dalang}, that 
	\begin{equation}\label{eq:hawkes}
		\lim_{|\xi|\to\infty}\frac{\Re\Psi(\xi)}{\log|\xi|}
		=\infty. 
	\end{equation}
	Then, for each fixed $t,\alpha>0$, the process $S_\alpha(t\,,\cdot)$
	has a modification that is in $C^\infty(\R)$. 
\end{proposition}

\begin{proof}
	By Plancherel's theorem, if $\mu$ is a finite signed measure on $\R$, then
	\begin{equation}
		\E\left(|S_\alpha(t\,,\mu)|^2\right) =
		\int_t^\infty\e^{-\alpha s}
		\| P^*_s\mu\|^2\,\d s
		= \frac{1}{2\pi}\int_t^\infty \e^{-\alpha s}\left\|
		\widehat{P^*_s\mu}\right\|^2\,\d s,
	\end{equation}
	where `` $\hat{}$ '' denotes the Fourier transform, normalized so that
	\begin{equation}
		\hat{g}(\xi) = \int_{-\infty}^\infty \e^{i\xi z}g(z)\,\d z
		\quad\text{for every $g\in L^1(\R$)}.
	\end{equation}
	But $\{P^*_s\}_{s\ge 0}$ is a convolution semigroup in the present
	setting, and has Fourier multiplier $\exp(-s\Psi(-\xi))$. That is,
	\begin{equation}
		\widehat{P^*_s\mu}(\xi)  = \e^{-s\Psi(-\xi)}\hat\mu(\xi)
		\quad\text{for all $\xi\in\R$ and $s\ge 0$}.
	\end{equation}
	Therefore, $\|\widehat{P^*_s\mu}\|^2 =
	\int_{-\infty}^\infty \e^{-2s\Re\Psi(\xi)}|\hat\mu(\xi)|^2\,\d\xi$.
	From this and the Tonelli theorem we deduce the following:
	\begin{equation}\label{eq:moments}
		\E\left(|S_\alpha(t\,,\mu)|^2\right) =
		\frac{1}{2\pi}\int_{-\infty}^\infty \frac{|\hat\mu(\xi)|^2}{
		\alpha+2\Re\Psi(\xi)}\e^{-(\alpha+2\Re\Psi(\xi))t}\,\d\xi.
	\end{equation}
	Recall \cite{GV} that the generalized $n$th derivative of $S_\alpha(t\,,\cdot)$
	is defined as the following random field:
	\begin{equation}
		S_\alpha^{(n)}(t\,,\phi) := (-1)^n S_\alpha(t\,,\phi^{(n)}),
	\end{equation}
	for all rapidly-decreasing test functions $\phi$ on $\R$. Here,
	$\phi^{(n)}$ denotes the $n$th derivative of $\phi$.
	Since the Fourier transform of $\phi^{(n)}$ is $i^n\xi^n\hat\phi(\xi)$,
	\eqref{eq:moments} implies that the following holds for all rapidly-decreasing
	test functions $\phi$ on $\R$:
	\begin{equation}\begin{split}
		\E\left(\left| S_\alpha^{(n)}(t\,,\phi) \right|^2\right)
			&\le \frac{1}{2\pi}\int_{-\infty}^\infty\frac{|\xi|^n\cdot |\hat\phi(\xi)|^2}{
			\alpha+2\Re\Psi(\xi)}\e^{-(\alpha+2\Re\Psi(\xi))t}\,\d\xi\\
		&\le\frac{\|\phi\|_{L^1(\R)}^2}{2\pi\alpha}\int_{-\infty}^\infty
			|\xi|^n\e^{-2t\Re\Psi(\xi)}\,\d\xi.
	\end{split}\end{equation}
	The growth condition on $\Re\Psi$ ensures that the integral is finite
	regardless of the value of $n$ and $t$. Therefore, a density argument
	shows that $Z^{(n)}_\alpha(t\,,\mu)$ can be defined as a $L^2(\P)$-continuous
	Gaussian random field indexed by all finite signed Borel measures $\mu$ on $\R$,
	and 
	\begin{equation}
		\E\left(\left| S_\alpha^{(n)}(t\,,\mu) \right|^2\right)
		\le\frac{(|\mu|(\R))^2}{2\pi\alpha}\int_{-\infty}^\infty
		|\xi|^{2n}\e^{-2t\Re\Psi(\xi)}\,\d\xi.
	\end{equation}
	This estimate and the Kolmogorov continuity theorem together imply
	that $\R\ni x\mapsto S_\alpha^{(n)}(t\,,x):=S_\alpha(t\,,\delta_x)$
	has a modification that is a bona fide continuous Gaussian process such that
	$S_\alpha^{(n)}(t\,,\mu)=\int S_\alpha^{(n)}(t\,,x)\,\mu(\d x)$. Apply this
	with $\mu$ replaced by a rapidly-decreasing test function $\phi$
	and apply integration by parts to deduce that $S_\alpha$ has a $C^\infty$
	modification.
\end{proof}

It would be interesting to know when \eqref{eq:dalang} 
implies \eqref{eq:hawkes}. This turns out to be a perplexing
problem, about which we next say a few words.

\begin{remark}\begin{enumerate}
	\item It is easy to see that Condition \eqref{eq:dalang} always implies that 		$\limsup_{|\xi|\to\infty}\Re\Psi(\xi)/|\xi|=\infty$. Therefore,
		in order to understand when \eqref{eq:dalang} implies \eqref{eq:hawkes}, we 
		need to know the $\liminf$ behavior of $\Re\Psi(\xi)/\log|\xi|$ for large
		values of $|\xi|$.
	\item It is shown in \cite[Lemma 8.1]{FKN} that \eqref{eq:dalang} implies the existence
		of transition densities $\{p_t(x)\}_{t>0,x\in\R}$ for $X$ such that
		$p_t$ is a bounded function for each fixed $t>0$. Theorem 6 of Hawkes \cite{Hawkes}
		then tell us that
		\begin{equation}\label{eq:hawkes:1}
			\lim_{|\xi|\to\infty}\frac{H(\xi)}{\log|\xi|}=\infty,
		\end{equation}
		where $H$ denotes Hardy--Littlewood's monotone increasing rearrangement
		of $\Re\Psi$. Condition \eqref{eq:hawkes:1} is tantalizingly close to
		\eqref{eq:hawkes}.
	\item It is easy to show that \eqref{eq:dalang} implies \eqref{eq:hawkes}
		if $\Re\Psi$ is nondecreasing on $[q\,,\infty)$ for $q$ sufficiently large.
		More generally, suppose $\Re\Psi$ is ``quasi-increasing'' on $\R_+$
		in the sense that there exists $C,q>0$ such that
		\begin{equation}
			\Re\Psi(2z)\ge C
			\sup_{u\in[z,2z]}\Re\Psi(u)\qquad
			\text{for all $z>q$}.
		\end{equation}
		Then \eqref{eq:dalang} implies \eqref{eq:hawkes}. Here is the [abelian]
		proof: For all $z>q$,
		\begin{equation}
			\int_z^{2z}\frac{\d\xi}{1+2\Re\Psi(\xi)}
			=2\int_{z/2}^z\frac{\d\xi}{1+2\Re\Psi(2\xi)}
			\ge \frac{z}{1+2C\Re\Psi(z)}
		\end{equation}
		Let $z$ tend to infinity and apply symmetry to find that
		\begin{equation}\label{eq:hawkes:bis}
			\lim_{|z|\to\infty}\frac{\Re\Psi(z)}{|z|}=\infty.
		\end{equation}
		Clearly, this [more than] implies \eqref{eq:hawkes}.
		\qed
\end{enumerate}\end{remark}

The preceding remarks suggest that, quite frequently, \eqref{eq:dalang}
implies \eqref{eq:hawkes}. We do not know whether or not \eqref{eq:dalang}
always implies \eqref{eq:hawkes}. But we are aware of some easy-to-check
conditions that guarantee this property. Let us state two such conditions next.
First, we recall the following two of the three well-known functions of Feller:
\begin{equation}
	K(\epsilon):=\epsilon^{-2}\int_{|z|\le\epsilon} z^2\,\nu(\d z)
	\quad\text{and}
	\quad G(\epsilon):=\nu\left\{x\in\R:\, |x|>\epsilon\right\},
\end{equation}
defined for all $\epsilon>0$. Then we have the following:

\begin{lemma}
	Suppose that at least one of the following two conditions holds:
	(i) $X$ has a nontrivial gaussian component; or (ii)
	\begin{equation}\label{eq:K/G}
		\limsup_{\epsilon\to 0^+}\frac{G(\epsilon)}{K(\epsilon)}<\infty.
	\end{equation}
	Then \eqref{eq:dalang} implies \eqref{eq:hawkes:bis}, 
	whence \eqref{eq:hawkes}.
\end{lemma}

Conditions \eqref{eq:K/G} and \eqref{eq:dalang} are not contradictory.
For a simple example, one can consider $X$ to be a symmetric stable
process of index $\alpha\in(0\,,2]$. Then, \eqref{eq:dalang}
holds if and only if $\alpha>1$. And, when $\alpha\in(1\,,2)$,
\eqref{eq:K/G} holds automatically;
in fact, in this case,
$G(\epsilon)$ and $K(\epsilon)$ both grow within constant
multiples of $\epsilon^{-\alpha}$
as $\epsilon\to 0^+$.

\begin{proof}
	In the case that $X$ contains a nontrivial gaussian component, the L\'evy--Khintchine
	formula \cite[p.\ 13]{Bertoin} implies that there exists $\sigma>0$ such that
	$\Re\Psi(\xi) \ge \sigma^2\xi^2$ for all $\xi\in\R$.
	Therefore, the result holds in this case. Thus, let us consider the case
	where there is no gaussian component in $X$ and \eqref{eq:K/G} holds.
	
	Because $1-\cos\theta\ge\theta^2/3$ for $\theta\in(0\,,1)$,
	we may apply the L\'evy--Khintchine formula to find the following
	well-known bound: For all $\xi>0$,
	\begin{equation}\label{eq:R1}
		\Re\Psi(\xi)=\int_{-\infty}^\infty
		(1-\cos(|x|\xi))\,\nu(\d\xi)\ge\frac13 K(1/\xi).
	\end{equation}
	
	Now we apply an averaging argument from harmonic analysis;
	define
	\begin{equation}
		\mathcal{R}(\xi) := \frac{1}{\xi}\int_0^\xi \Re\Psi(z)\,\d z
		\quad\text{for all $\xi>0$}.
	\end{equation}
	By the L\'evy--Khintchine formula,
	\begin{equation}
		\mathcal{R}(\xi) = \int_{-\infty}^\infty( 1-\mathop{\text{sinc}}(|x|\xi)\,\nu(\d x),
	\end{equation}
	where $\mathop{\text{sinc}}\theta:=\sin\theta/\theta$, as usual.
	Because $1-\mathop{\text{sinc}}\theta\le \min(\theta^2/2\,,1)$
	for all $\theta>0$, it follows that for all $\xi>0$,
	\begin{equation}\label{eq:RKG}
		\mathcal{R}(\xi) \le \frac12 K(1/\xi) + G(1/\xi).
	\end{equation}
	Because of \eqref{eq:K/G}, \eqref{eq:R1}, and \eqref{eq:RKG},
	$\mathcal{R}(\xi) = O(\Re\Psi(\xi))$ as $\xi\to\infty$.
	Thanks to symmetry, it suffices to prove that
	\begin{equation}\label{eq:hawkes:bis1}
		\lim_{\xi\to\infty}\frac{\mathcal{R}(\xi)}{\xi}=\infty.
	\end{equation}
	To this end, we observe that for all $\alpha,\xi >0$,
	\begin{equation}\begin{split}
		\alpha+\mathcal{R}(\xi) &=\frac{1}{\xi}\int_0^\xi (\alpha+\Re\Psi(z))\,
			\d z\\
		&\ge\left(\frac1\xi\int_0^\xi \frac{\d z}{\alpha+2\Re\Psi(z)}\right)^{-1};
	\end{split}\end{equation}
	this follows from the Cauchy--Schwarz inequality. Consequently,
	\begin{equation}
		\liminf_{\xi\to\infty} \frac{\mathcal{R}(\xi)}{\xi}
		\ge \sup_{\alpha>0}\left(\int_0^\infty
		\frac{\d z}{\alpha+2\Re\Psi(z)}\right)^{-1}=\infty.
	\end{equation}
	Because of this and symmetry we obtain \eqref{eq:hawkes:bis1},
	and hence the lemma.
\end{proof}

\begin{remark}
	Our proof was based on the general idea that, under \eqref{eq:K/G},
	Condition \eqref{eq:hawkes}  implies \eqref{eq:hawkes:bis1}.
	The converse holds even without \eqref{eq:K/G}. In fact, because
	$1-\cos\theta\le \text{const}\cdot(1-\mathop{\text{sinc}}\theta)$
	for a universal constant, the L\'evy--Khintchine formula implies
	that $\Re\Psi(\xi)\le \text{const}\cdot\mathcal{R}(\xi)$ for the
	same universal constant.
	\qed
\end{remark}

We conclude this paper by presenting an example.

\begin{example}
	Consider the following stochastic cable equation: For $t>0$ and $x\in\R$,
	\begin{equation}\label{cable:beta}
		\frac{\partial}{\partial t} V_\alpha(t\,,x) = (\Delta_{\beta/2}
		V_\alpha)(t\,,x) - \frac{\alpha}{2}V_\alpha(t\,,x)+\dot{W}(t\,,x),
	\end{equation}
	with $V_\alpha(0\,,x)=0$. Here, $\Delta_{\beta/2}$ denotes the
	fractional Laplacian of index $\beta/2$, and $\beta\in(1\,,2]$, 
	normalized so that
	$\Delta_1f = f''$. Because $\beta>1$, \eqref{eq:dalang} is verified. 
	It is well known that $\Delta_{\beta/2}$ is the
	generator of $\{cX_t\}_{t\ge 0}$, where $X$ denote the symmetric
	stable process of index $\beta$ and $c=c(\beta)$ is a certain positive constant.
	In this case, $\Psi(\xi)=\text{const}\cdot|\xi|^\beta$ is real, and \eqref{eq:hawkes}
	holds. According to
	Proposition \ref{pr:smooth}, $S_\alpha(t\,,\cdot)\in C^\infty(\R)$
	[up to a modification] for every fixed $t>0$. Therefore, the solution
	$V_\alpha(t\,,\cdot)$ to the stochastic cable equation \eqref{cable} is a $C^\infty$
	perturbation of the Gaussian process $\eta_\alpha(t\,,\cdot)$
	that is associated to the local times of a symmetric stable process of index $\beta$.
	But $\eta_\alpha(t\,,\cdot)$ is fractional Brownian motion of index $\beta-1
	\in(0\,,1]$;
	see Marcus and Rosen \cite[p.\ 498]{MR}. Recall that $\eta_\alpha(t\,,\cdot)$ is a 
	continuous centered Gaussian process with
	\begin{equation}
		\E\left(|\eta_\alpha(t\,,x)-\eta_\alpha(t\,,y)|^2\right)=\text{const}
		\cdot |x-y|^{\beta-1}
		\quad\text{for $x,y\in\R$}.
	\end{equation}
	
	As a consequence, we find that local regularity of
	the solution to the cable equation \eqref{cable:beta} [in the space
	variable] is the same as local regularity of fractional Brownian motion
	of index $\beta-1$. 
	A similar phenomenon has
	been observed by Walsh \cite[Exercise 3.10, p.\ 326]{Walsh} in the case that we solve
	\eqref{cable:beta} for $\beta=2$ on a finite $x$-interval.
	Dynkin's isomorphism theorem \cite[Chapter 8]{MR} then implies that 
	the local regularity, in the space variable, of the solution to
	\eqref{cable:beta} is the same as the local regularity of the local time
	of a symmetric stable process of index $\beta$ in the space variable.
	And thanks to a change of measure (Proposition
	\ref{pr:girsanov}), the same is true of local regularity of the solution
	to the heat equation,
	\begin{equation}\left|\begin{split}
		&\frac{\partial}{\partial t}U(t\,,x) = (\Delta_{\beta/2}U)(t\,,x) + \dot{W}(t\,,x),\\
		&U(0\,,x)=0.
	\end{split}\right.\end{equation}
	
	One can consult other local-time examples of this general type by
	drawing upon the examples from the book of Marcus and Rosen \cite{MR}.
\qed
\end{example}

\noindent\textbf{Acknowledgements.} 
Many thanks are owed to several people and sources.
Professor Zhan Shi shared with us
his ideas generously. Among many other things,
Corollary \ref{co:bd2} was formulated with his help.
Professor Steven Evans generously shared with us his copy of
the unpublished manuscript by Hawkes \cite{Hawkes}. And
the University of Paris [13 and 6] kindly hosted D.K.
for a three-week visit, during which some of this research was
completed. We express our gratitude to them all.

\begin{spacing}{.8}
\begin{small}
\medskip
\begin{spacing}{.5}
\noindent\textbf{Nathalie Eisenbaum}
\noindent  Laboratoire de Probabilit\'es et Mod\`eles Al\'eatoires ,
	Universit\'e Paris VI - case 188  4, Place Jussieu, 75252,
	Paris Cedex 05, France.\\
\noindent\emph{Email:} \texttt{nathalie.eisenbaum@umpc.fr}\\[3mm]

\noindent\textbf{Mohammud Foondun}
\noindent School of Mathematics, Loughborough University,
		Leicestershire, LE11 3TU, UK
\noindent\emph{Email:} \texttt{M.I.Foondun@lboro.ac.uk}\\[3mm]
\noindent\textbf{Davar Khoshnevisan}
\noindent Department of Mathematics, University of Utah,
		Salt Lake City, UT 84112-0090
\noindent\emph{Email:} \texttt{davar@math.utah.edu}
\end{spacing}\end{small}
\end{spacing}

\end{document}